\documentclass{article}

\usepackage{tikz, amssymb, amsmath, amsthm, enumitem, skak, graphicx}
\usepackage[pdfborder = {0 0 0}]{hyperref}

\usepackage{multicol}
\usepackage{caption}
\usetikzlibrary{calc}
\usepackage{wrapfig}

\theoremstyle{definition}

\DeclareRobustCommand{\qbinom}{\genfrac{[}{]}{0pt}{}}
\DeclareMathOperator{\inv}{inv}
\DeclareMathOperator{\maj}{maj}
\DeclareMathOperator{\des}{des}

\newcommand{\patterna}{\!\!
  \begin{tikzpicture}[scale = .15, baseline = 1ex]
    \draw [color = black!25, thin] (0,0) grid (3,3);
    \draw [color=teal, fill=teal, thick] (0.5,1.5) circle (1.25ex);
    \draw [color=teal, fill=teal, thick] (1.5,2.5) circle (1.25ex);
    \draw [color=teal, fill=teal, thick] (2.5,0.5) circle (1.25ex);
    \draw [line width=0.01mm, thick, teal] (0,1) -- (3,1);
    \draw [line width=0.01mm, thick, teal] (1,0) -- (1,3);
  \end{tikzpicture}\,\,}

\definecolor{teal}{RGB}{102,153,204}
\definecolor{pink}{RGB}{204,102,119}

\title{The combinatorics of Jeff Remmel}
\author{Sergey Kitaev
  \footnote{Department of Mathematics and Statistics, University of
    Strathclyde, Glasgow, G1 1XH, UK; Email: sergey.kitaev@strath.ac.uk}
  \and
  Anthony Mendes
  \footnote{Department of Mathematics, California Polytechnic
    State University, San Luis Obispo, USA; Email: aamendes@calpoly.edu}}

\begin{document}

\maketitle

\begin{abstract}
  We give a brief overview of the life and combinatorics of Jeff Remmel, a
  mathematician with successful careers in both logic and combinatorics.
\end{abstract}

\section{Biography}

\begin{wrapfigure}{l}{0.33\textwidth}
\includegraphics[width=.9\linewidth]{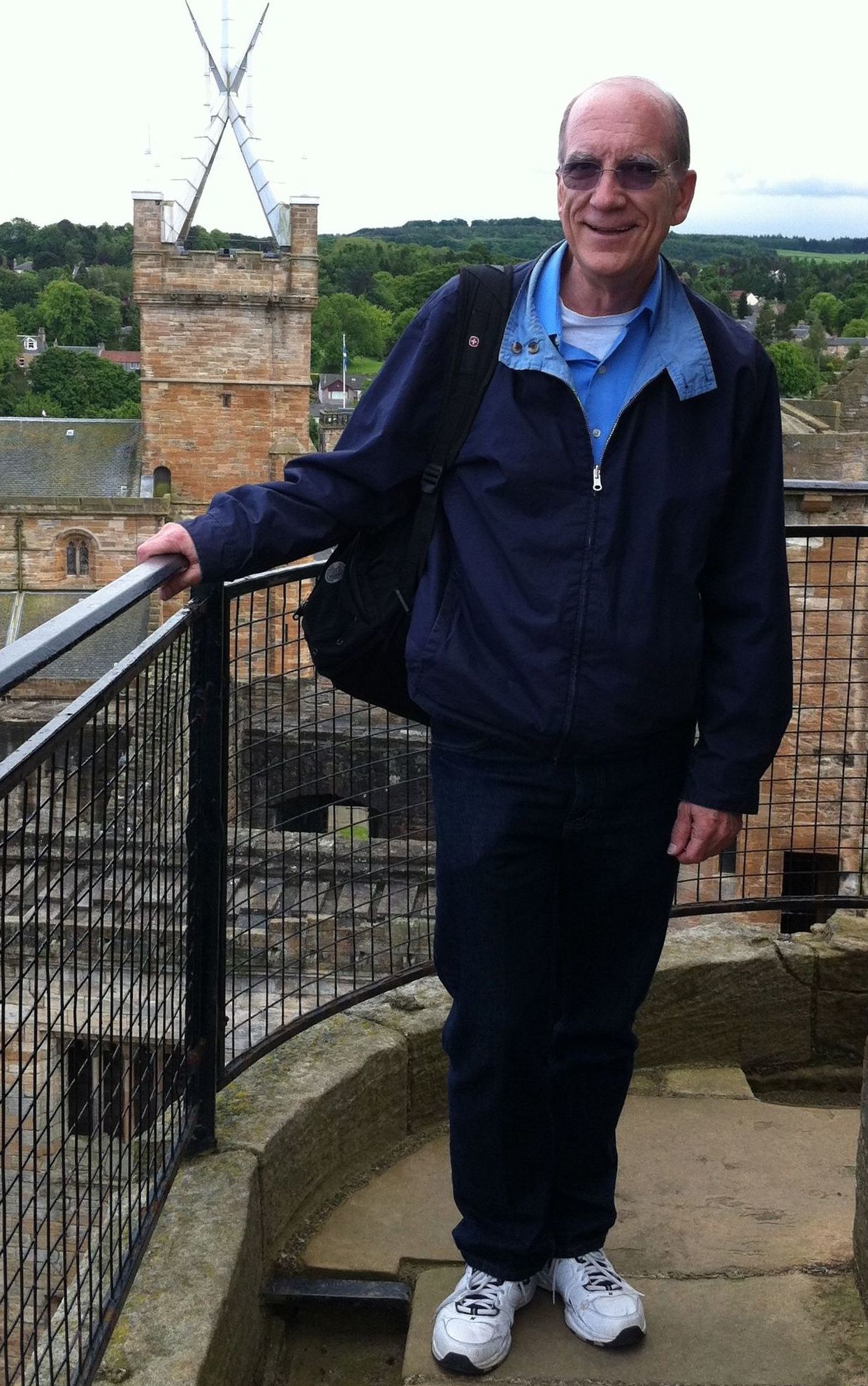}
\caption*{\begin{small}Linlithgow Palace, Scotland, June 2012.  Jeff visits the
    $10^{\mbox{\tiny{th}}}$ Permutation Patterns Conference hosted by the
    University of Strathclyde in Glasgow. The photo was taken by the first
    author.\end{small}}
\end{wrapfigure}

Born October 12, 1948 in Clintonville, Wisconsin, Jeff Remmel earned his
undergraduate degree from Swarthmore College in 1970 and his PhD in logic from
Cornell University in 1974.

Jeff was hired as an Assistant Professor in the Department of Mathematics at UC
San Diego at age 25, without officially finishing his PhD and without having
published a single paper!  He stumbled into the position after his thesis
advisor, Anil Nerode, recommended Jeff as a replacement hire for another
logician who rescinded the job offer late.  Jeff called his hiring a ``fluke
that will never happen again'' and that he was ``completely clueless'' about
the hiring and promotion process until years after working at the university.

Jeff's interview process was held completely over the phone and he first
visited the campus in 1974 when he arrived to teach his Fall courses.  He spent
the next 42 years as a faculty member at UC San Diego, building an
exceptionally successful academic career.

A perennial favorite among mathematics graduate students, Jeff enjoyed teaching
the introductory graduate courses in enumerative and algebraic combinatorics.
He frequently volunteered to teach courses beyond his assigned course load.
This provided him with a constant stream of PhD students.  Jeff took students
freely and without reservation, immediately accepting anyone who asked to work
with him.  He graduated 33 PhD students with an unusually high number of his
students becoming university faculty members themselves.

Jeff would playfully joke about the frustrations of dealing with
administrators.  During a lecture on derangements in 2001, he said ``The old
story goes, you have a dumb blonde as a hat check girl---it could be a brunette,
or it could be a male.  If it were an administrator, they'd really be dumb'',
following with ``I will deny the remark about administrators if it ever leaves
the room''.  While he was department chair, Jeff told departing graduate
students ``Don't become department chair'' with his tongue-in-cheek.

These jokes were just a facade. Jeff truly relished his many leadership roles
at the university.  He was department chair for 4 years, associate dean for 16
years, and interim dean for one year.  He additionally was a founding director
for a state-wide program to train future K--12 teachers (CalTeach), was a
founding director for a program on improving undergraduate education, was a
director of a Summer Bridge program to the university, helped create a data
science major, and was instrumental in hiring founding faculty for an MBA
program.

Despite the time and energy spent on teaching and administrative work, Jeff was
an amazingly prolific mathematician, publishing 322 research articles with over
100 coauthors.  This chart contains one square for each publication:

\begin{center}
  \begin{tikzpicture}[scale = .25]
    \node at (76.5,-1)  {\scriptsize 1976};
    \node at (80.5,-1)  {\scriptsize 1980};
    \node at (85.5,-1)  {\scriptsize 1985};
    \node at (90.5,-1)  {\scriptsize 1990};
    \node at (95.5,-1)  {\scriptsize 1995};
    \node at (100.5,-1) {\scriptsize 2000};
    \node at (105.5,-1) {\scriptsize 2005};
    \node at (110.5,-1) {\scriptsize 2010};
    \node at (115.5,-1) {\scriptsize 2015};
    \node at (120.5,-1) {\scriptsize 2020};
    \fill [color = pink] (76,0) rectangle ++(1,1);
    \fill [color = pink] (76,1) rectangle ++(1,1);
    \fill [color = pink] (77,0) rectangle ++(1,1);
    \fill [color = pink] (78,0) rectangle ++(1,1);
    \fill [color = pink] (78,1) rectangle ++(1,1);
    \fill [color = pink] (78,2) rectangle ++(1,1);
    \fill [color = pink] (79,0) rectangle ++(1,1);
    \fill [color = pink] (79,1) rectangle ++(1,1);
    \fill [color = teal] (79,2) rectangle ++(1,1);
    \fill [color = pink] (80,0) rectangle ++(1,1);
    \fill [color = pink] (80,1) rectangle ++(1,1);
    \fill [color = pink] (80,2) rectangle ++(1,1);
    \fill [color = pink] (80,3) rectangle ++(1,1);
    \fill [color = pink] (80,4) rectangle ++(1,1);
    \fill [color = teal] (80,5) rectangle ++(1,1);
    \fill [color = pink] (81,0) rectangle ++(1,1);
    \fill [color = pink] (81,1) rectangle ++(1,1);
    \fill [color = pink] (81,2) rectangle ++(1,1);
    \fill [color = pink] (81,3) rectangle ++(1,1);
    \fill [color = pink] (81,4) rectangle ++(1,1);
    \fill [color = pink] (81,5) rectangle ++(1,1);
    \fill [color = pink] (81,6) rectangle ++(1,1);
    \fill [color = pink] (81,7) rectangle ++(1,1);
    \fill [color = teal] (81,8) rectangle ++(1,1);
    \fill [color = pink] (82,0) rectangle ++(1,1);
    \fill [color = pink] (82,1) rectangle ++(1,1);
    \fill [color = pink] (82,2) rectangle ++(1,1);
    \fill [color = teal] (82,3) rectangle ++(1,1);
    \fill [color = teal] (82,4) rectangle ++(1,1);
    \fill [color = pink] (83,0) rectangle ++(1,1);
    \fill [color = pink] (83,1) rectangle ++(1,1);
    \fill [color = pink] (83,2) rectangle ++(1,1);
    \fill [color = pink] (83,3) rectangle ++(1,1);
    \fill [color = teal] (83,4) rectangle ++(1,1);
    \fill [color = teal] (83,5) rectangle ++(1,1);
    \fill [color = pink] (84,0) rectangle ++(1,1);
    \fill [color = pink] (84,1) rectangle ++(1,1);
    \fill [color = pink] (84,2) rectangle ++(1,1);
    \fill [color = pink] (84,3) rectangle ++(1,1);
    \fill [color = teal] (84,4) rectangle ++(1,1);
    \fill [color = teal] (84,5) rectangle ++(1,1);
    \fill [color = teal] (84,6) rectangle ++(1,1);
    \fill [color = teal] (84,7) rectangle ++(1,1);
    \fill [color = pink] (85,0) rectangle ++(1,1);
    \fill [color = pink] (85,1) rectangle ++(1,1);
    \fill [color = pink] (85,2) rectangle ++(1,1);
    \fill [color = teal] (85,3) rectangle ++(1,1);
    \fill [color = teal] (85,4) rectangle ++(1,1);
    \fill [color = teal] (85,5) rectangle ++(1,1);
    \fill [color = teal] (85,6) rectangle ++(1,1);
    \fill [color = pink] (86,0) rectangle ++(1,1);
    \fill [color = pink] (86,1) rectangle ++(1,1);
    \fill [color = teal] (86,2) rectangle ++(1,1);
    \fill [color = teal] (86,3) rectangle ++(1,1);
    \fill [color = teal] (86,4) rectangle ++(1,1);
    \fill [color = teal] (86,5) rectangle ++(1,1);
    \fill [color = pink] (87,0) rectangle ++(1,1);
    \fill [color = pink] (87,1) rectangle ++(1,1);
    \fill [color = pink] (87,2) rectangle ++(1,1);
    \fill [color = pink] (87,3) rectangle ++(1,1);
    \fill [color = pink] (87,4) rectangle ++(1,1);
    \fill [color = teal] (87,5) rectangle ++(1,1);
    \fill [color = teal] (87,6) rectangle ++(1,1);
    \fill [color = teal] (88,0) rectangle ++(1,1);
    \fill [color = teal] (88,1) rectangle ++(1,1);
    \fill [color = pink] (89,0) rectangle ++(1,1);
    \fill [color = pink] (89,1) rectangle ++(1,1);
    \fill [color = pink] (89,2) rectangle ++(1,1);
    \fill [color = pink] (89,3) rectangle ++(1,1);
    \fill [color = teal] (89,4) rectangle ++(1,1);
    \fill [color = pink] (90,0) rectangle ++(1,1);
    \fill [color = pink] (90,1) rectangle ++(1,1);
    \fill [color = pink] (90,2) rectangle ++(1,1);
    \fill [color = pink] (90,3) rectangle ++(1,1);
    \fill [color = pink] (90,4) rectangle ++(1,1);
    \fill [color = pink] (90,5) rectangle ++(1,1);
    \fill [color = teal] (90,6) rectangle ++(1,1);
    \fill [color = teal] (90,7) rectangle ++(1,1);
    \fill [color = teal] (90,8) rectangle ++(1,1);
    \fill [color = teal] (90,9) rectangle ++(1,1);
    \fill [color = pink] (91,0) rectangle ++(1,1);
    \fill [color = pink] (91,1) rectangle ++(1,1);
    \fill [color = teal] (91,2) rectangle ++(1,1);
    \fill [color = teal] (91,3) rectangle ++(1,1);
    \fill [color = pink] (92,0) rectangle ++(1,1);
    \fill [color = pink] (92,1) rectangle ++(1,1);
    \fill [color = pink] (92,2) rectangle ++(1,1);
    \fill [color = pink] (92,3) rectangle ++(1,1);
    \fill [color = pink] (92,4) rectangle ++(1,1);
    \fill [color = pink] (92,5) rectangle ++(1,1);
    \fill [color = teal] (92,6) rectangle ++(1,1);
    \fill [color = teal] (92,7) rectangle ++(1,1);
    \fill [color = pink] (93,0) rectangle ++(1,1);
    \fill [color = pink] (93,1) rectangle ++(1,1);
    \fill [color = pink] (93,2) rectangle ++(1,1);
    \fill [color = pink] (93,3) rectangle ++(1,1);
    \fill [color = teal] (93,4) rectangle ++(1,1);
    \fill [color = teal] (93,5) rectangle ++(1,1);
    \fill [color = pink] (94,0) rectangle ++(1,1);
    \fill [color = pink] (94,1) rectangle ++(1,1);
    \fill [color = pink] (94,2) rectangle ++(1,1);
    \fill [color = pink] (94,3) rectangle ++(1,1);
    \fill [color = pink] (94,4) rectangle ++(1,1);
    \fill [color = teal] (94,5) rectangle ++(1,1);
    \fill [color = teal] (94,6) rectangle ++(1,1);
    \fill [color = pink] (95,0) rectangle ++(1,1);
    \fill [color = pink] (95,1) rectangle ++(1,1);
    \fill [color = pink] (95,2) rectangle ++(1,1);
    \fill [color = pink] (95,3) rectangle ++(1,1);
    \fill [color = pink] (95,4) rectangle ++(1,1);
    \fill [color = pink] (95,5) rectangle ++(1,1);
    \fill [color = pink] (95,6) rectangle ++(1,1);
    \fill [color = pink] (95,7) rectangle ++(1,1);
    \fill [color = pink] (95,8) rectangle ++(1,1);
    \fill [color = pink] (95,9) rectangle ++(1,1);
    \fill [color = black!66] (95,10) rectangle ++(1,1);
    \fill [color = teal] (95,11) rectangle ++(1,1);
    \fill [color = teal] (95,12) rectangle ++(1,1);
    \fill [color = teal] (95,13) rectangle ++(1,1);
    \fill [color = pink] (96,0) rectangle ++(1,1);
    \fill [color = pink] (96,1) rectangle ++(1,1);
    \fill [color = pink] (96,2) rectangle ++(1,1);
    \fill [color = pink] (96,3) rectangle ++(1,1);
    \fill [color = pink] (96,4) rectangle ++(1,1);
    \fill [color = pink] (96,5) rectangle ++(1,1);
    \fill [color = pink] (96,6) rectangle ++(1,1);
    \fill [color = pink] (96,7) rectangle ++(1,1);
    \fill [color = teal] (96,8) rectangle ++(1,1);
    \fill [color = teal] (96,9) rectangle ++(1,1);
    \fill [color = teal] (96,10) rectangle ++(1,1);
    \fill [color = pink] (97,0) rectangle ++(1,1);
    \fill [color = pink] (97,1) rectangle ++(1,1);
    \fill [color = pink] (97,2) rectangle ++(1,1);
    \fill [color = pink] (97,3) rectangle ++(1,1);
    \fill [color = pink] (97,4) rectangle ++(1,1);
    \fill [color = pink] (97,5) rectangle ++(1,1);
    \fill [color = pink] (97,6) rectangle ++(1,1);
    \fill [color = pink] (97,7) rectangle ++(1,1);
    \fill [color = teal] (97,8) rectangle ++(1,1);
    \fill [color = pink] (98,0) rectangle ++(1,1);
    \fill [color = pink] (98,1) rectangle ++(1,1);
    \fill [color = pink] (98,2) rectangle ++(1,1);
    \fill [color = pink] (98,3) rectangle ++(1,1);
    \fill [color = pink] (98,4) rectangle ++(1,1);
    \fill [color = pink] (98,5) rectangle ++(1,1);
    \fill [color = teal] (98,6) rectangle ++(1,1);
    \fill [color = teal] (98,7) rectangle ++(1,1);
    \fill [color = teal] (98,8) rectangle ++(1,1);
    \fill [color = teal] (98,9) rectangle ++(1,1);
    \fill [color = pink] (99,0) rectangle ++(1,1);
    \fill [color = pink] (99,1) rectangle ++(1,1);
    \fill [color = pink] (99,2) rectangle ++(1,1);
    \fill [color = pink] (99,3) rectangle ++(1,1);
    \fill [color = pink] (99,4) rectangle ++(1,1);
    \fill [color = pink] (99,5) rectangle ++(1,1);
    \fill [color = teal] (99,6) rectangle ++(1,1);
    \fill [color = pink] (100,0) rectangle ++(1,1);
    \fill [color = pink] (101,0) rectangle ++(1,1);
    \fill [color = teal] (101,1) rectangle ++(1,1);
    \fill [color = pink] (102,0) rectangle ++(1,1);
    \fill [color = pink] (102,1) rectangle ++(1,1);
    \fill [color = pink] (102,2) rectangle ++(1,1);
    \fill [color = teal] (102,3) rectangle ++(1,1);
    \fill [color = pink] (103,0) rectangle ++(1,1);
    \fill [color = pink] (103,1) rectangle ++(1,1);
    \fill [color = pink] (103,2) rectangle ++(1,1);
    \fill [color = black!66] (103,3) rectangle ++(1,1);
    \fill [color = teal] (103,4) rectangle ++(1,1);
    \fill [color = teal] (103,5) rectangle ++(1,1);
    \fill [color = pink] (104,0) rectangle ++(1,1);
    \fill [color = pink] (104,1) rectangle ++(1,1);
    \fill [color = pink] (104,2) rectangle ++(1,1);
    \fill [color = pink] (104,3) rectangle ++(1,1);
    \fill [color = pink] (104,4) rectangle ++(1,1);
    \fill [color = pink] (104,5) rectangle ++(1,1);
    \fill [color = teal] (104,6) rectangle ++(1,1);
    \fill [color = teal] (104,7) rectangle ++(1,1);
    \fill [color = teal] (104,8) rectangle ++(1,1);
    \fill [color = teal] (104,9) rectangle ++(1,1);
    \fill [color = teal] (104,10) rectangle ++(1,1);
    \fill [color = pink] (105,0) rectangle ++(1,1);
    \fill [color = pink] (105,1) rectangle ++(1,1);
    \fill [color = pink] (105,2) rectangle ++(1,1);
    \fill [color = pink] (105,3) rectangle ++(1,1);
    \fill [color = teal] (105,4) rectangle ++(1,1);
    \fill [color = teal] (105,5) rectangle ++(1,1);
    \fill [color = pink] (106,0) rectangle ++(1,1);
    \fill [color = pink] (106,1) rectangle ++(1,1);
    \fill [color = pink] (106,2) rectangle ++(1,1);
    \fill [color = teal] (106,3) rectangle ++(1,1);
    \fill [color = teal] (106,4) rectangle ++(1,1);
    \fill [color = teal] (106,5) rectangle ++(1,1);
    \fill [color = teal] (106,6) rectangle ++(1,1);
    \fill [color = teal] (106,7) rectangle ++(1,1);
    \fill [color = teal] (106,8) rectangle ++(1,1);
    \fill [color = teal] (106,9) rectangle ++(1,1);
    \fill [color = pink] (107,0) rectangle ++(1,1);
    \fill [color = pink] (107,1) rectangle ++(1,1);
    \fill [color = pink] (107,2) rectangle ++(1,1);
    \fill [color = teal] (107,3) rectangle ++(1,1);
    \fill [color = pink] (108,0) rectangle ++(1,1);
    \fill [color = pink] (108,1) rectangle ++(1,1);
    \fill [color = pink] (108,2) rectangle ++(1,1);
    \fill [color = pink] (108,3) rectangle ++(1,1);
    \fill [color = pink] (108,4) rectangle ++(1,1);
    \fill [color = pink] (108,5) rectangle ++(1,1);
    \fill [color = black!66] (108,6) rectangle ++(1,1);
    \fill [color = teal] (108,7) rectangle ++(1,1);
    \fill [color = teal] (108,8) rectangle ++(1,1);
    \fill [color = teal] (108,9) rectangle ++(1,1);
    \fill [color = teal] (108,10) rectangle ++(1,1);
    \fill [color = teal] (108,11) rectangle ++(1,1);
    \fill [color = teal] (108,12) rectangle ++(1,1);
    \fill [color = pink] (109,0) rectangle ++(1,1);
    \fill [color = pink] (109,1) rectangle ++(1,1);
    \fill [color = pink] (109,2) rectangle ++(1,1);
    \fill [color = pink] (109,3) rectangle ++(1,1);
    \fill [color = pink] (109,4) rectangle ++(1,1);
    \fill [color = pink] (109,5) rectangle ++(1,1);
    \fill [color = pink] (109,6) rectangle ++(1,1);
    \fill [color = teal] (109,7) rectangle ++(1,1);
    \fill [color = teal] (109,8) rectangle ++(1,1);
    \fill [color = teal] (109,9) rectangle ++(1,1);
    \fill [color = teal] (109,10) rectangle ++(1,1);
    \fill [color = teal] (109,11) rectangle ++(1,1);
    \fill [color = pink] (110,0) rectangle ++(1,1);
    \fill [color = teal] (110,1) rectangle ++(1,1);
    \fill [color = teal] (110,2) rectangle ++(1,1);
    \fill [color = teal] (110,3) rectangle ++(1,1);
    \fill [color = teal] (110,4) rectangle ++(1,1);
    \fill [color = teal] (110,5) rectangle ++(1,1);
    \fill [color = teal] (110,6) rectangle ++(1,1);
    \fill [color = teal] (110,7) rectangle ++(1,1);
    \fill [color = teal] (110,8) rectangle ++(1,1);
    \fill [color = pink] (111,0) rectangle ++(1,1);
    \fill [color = pink] (111,1) rectangle ++(1,1);
    \fill [color = pink] (111,2) rectangle ++(1,1);
    \fill [color = pink] (111,3) rectangle ++(1,1);
    \fill [color = pink] (111,4) rectangle ++(1,1);
    \fill [color = pink] (111,5) rectangle ++(1,1);
    \fill [color = pink] (111,6) rectangle ++(1,1);
    \fill [color = teal] (111,7) rectangle ++(1,1);
    \fill [color = teal] (111,8) rectangle ++(1,1);
    \fill [color = teal] (111,9) rectangle ++(1,1);
    \fill [color = teal] (111,10) rectangle ++(1,1);
    \fill [color = teal] (111,11) rectangle ++(1,1);
    \fill [color = teal] (111,12) rectangle ++(1,1);
    \fill [color = teal] (111,13) rectangle ++(1,1);
    \fill [color = teal] (111,14) rectangle ++(1,1);
    \fill [color = teal] (111,15) rectangle ++(1,1);
    \fill [color = teal] (111,16) rectangle ++(1,1);
    \fill [color = pink] (112,0) rectangle ++(1,1);
    \fill [color = pink] (112,1) rectangle ++(1,1);
    \fill [color = pink] (112,2) rectangle ++(1,1);
    \fill [color = pink] (112,3) rectangle ++(1,1);
    \fill [color = teal] (112,4) rectangle ++(1,1);
    \fill [color = teal] (112,5) rectangle ++(1,1);
    \fill [color = teal] (112,6) rectangle ++(1,1);
    \fill [color = teal] (112,7) rectangle ++(1,1);
    \fill [color = teal] (112,8) rectangle ++(1,1);
    \fill [color = teal] (112,9) rectangle ++(1,1);
    \fill [color = teal] (112,10) rectangle ++(1,1);
    \fill [color = teal] (112,11) rectangle ++(1,1);
    \fill [color = pink] (113,0) rectangle ++(1,1);
    \fill [color = pink] (113,1) rectangle ++(1,1);
    \fill [color = pink] (113,2) rectangle ++(1,1);
    \fill [color = black!66] (113,3) rectangle ++(1,1);
    \fill [color = teal] (113,4) rectangle ++(1,1);
    \fill [color = teal] (113,5) rectangle ++(1,1);
    \fill [color = teal] (113,6) rectangle ++(1,1);
    \fill [color = teal] (113,7) rectangle ++(1,1);
    \fill [color = teal] (113,8) rectangle ++(1,1);
    \fill [color = teal] (113,9) rectangle ++(1,1);
    \fill [color = teal] (113,10) rectangle ++(1,1);
    \fill [color = pink] (114,0) rectangle ++(1,1);
    \fill [color = pink] (114,0) rectangle ++(1,1);
    \fill [color = pink] (114,1) rectangle ++(1,1);
    \fill [color = teal] (114,2) rectangle ++(1,1);
    \fill [color = teal] (114,3) rectangle ++(1,1);
    \fill [color = teal] (114,4) rectangle ++(1,1);
    \fill [color = teal] (114,5) rectangle ++(1,1);
    \fill [color = teal] (114,6) rectangle ++(1,1);
    \fill [color = pink] (115,0) rectangle ++(1,1);
    \fill [color = teal] (115,1) rectangle ++(1,1);
    \fill [color = teal] (115,2) rectangle ++(1,1);
    \fill [color = teal] (115,3) rectangle ++(1,1);
    \fill [color = teal] (115,4) rectangle ++(1,1);
    \fill [color = teal] (115,5) rectangle ++(1,1);
    \fill [color = teal] (115,6) rectangle ++(1,1);
    \fill [color = teal] (115,7) rectangle ++(1,1);
    \fill [color = teal] (115,8) rectangle ++(1,1);
    \fill [color = teal] (115,9) rectangle ++(1,1);
    \fill [color = teal] (115,10) rectangle ++(1,1);
    \fill [color = pink] (116,0) rectangle ++(1,1);
    \fill [color = teal] (116,1) rectangle ++(1,1);
    \fill [color = teal] (116,2) rectangle ++(1,1);
    \fill [color = teal] (116,3) rectangle ++(1,1);
    \fill [color = teal] (116,4) rectangle ++(1,1);
    \fill [color = teal] (116,5) rectangle ++(1,1);
    \fill [color = teal] (116,6) rectangle ++(1,1);
    \fill [color = pink] (117,0) rectangle ++(1,1);
    \fill [color = black!66] (117,1) rectangle ++(1,1);
    \fill [color = teal] (117,2) rectangle ++(1,1);
    \fill [color = teal] (117,3) rectangle ++(1,1);
    \fill [color = teal] (117,4) rectangle ++(1,1);
    \fill [color = teal] (117,5) rectangle ++(1,1);
    \fill [color = teal] (117,6) rectangle ++(1,1);
    \fill [color = teal] (118,0) rectangle ++(1,1);
    \fill [color = teal] (118,1) rectangle ++(1,1);
    \fill [color = teal] (118,2) rectangle ++(1,1);
    \fill [color = teal] (119,0) rectangle ++(1,1);
    \fill [color = teal] (119,1) rectangle ++(1,1);
    \fill [color = teal] (119,2) rectangle ++(1,1);
    \fill [color = teal] (119,3) rectangle ++(1,1);
    \fill [color = teal] (119,4) rectangle ++(1,1);
    \fill [color = teal] (119,5) rectangle ++(1,1);
    \fill [color = teal] (119,6) rectangle ++(1,1);
    \fill [color = teal] (119,7) rectangle ++(1,1);
    \fill [color = pink] (120,0) rectangle ++(1,1);
    \fill [color = pink] (120,1) rectangle ++(1,1);
    \fill [color = teal] (120,2) rectangle ++(1,1);
    \draw [color = white, ultra thin] (76,0) grid (121,17);
    \draw [color = black!15, ultra thin] (76,0) -- (121,0);
    \draw [color = black!15, ultra thin] (76,5) -- (121,5);
    \draw [color = black!15, ultra thin] (76,10) -- (121,10);
    \draw [color = black!15, ultra thin] (76,15) -- (121,15);
    \node at (75,5) {\scriptsize 5};
    \node at (75,10) {\scriptsize 10};
    \node at (75,15) {\scriptsize 15};
  \end{tikzpicture}
\end{center}
The 167 pink squares indicate publications in logic, the 149 teal squares (the
color of the poplin shirt he often wore to work) indicate publications in
combinatorics, and the 5 gray squares indicate publications that are neither
logic or combinatorics.  His publications would tend to be quite lengthy, with
many papers over 30 pages.

Two research papers per year is considered a very respectable rate of
publication for a research mathematician, giving an expected 88 publications
over the span of 44 years \cite{AMS}.  This means that Jeff had two almost
completely separate remarkably productive careers; one in logic and one in
combinatorics.

This paper highlights Jeff's accomplishments in combinatorics only, with the
remaining sections outlining Jeff's combinatorics results sorted by theme.
Jeff discovered algebraic and enumerative combinatorics after taking a course
from his longtime colleague Adriano Garsia while Jeff was an assistant
professor.  Even though Jeff was indeed a strong researcher in logic, he
probably ended up more well known for combinatorics with the great majority of
his PhD students, invited talks and grants on the subject.

Personally, Jeff was generous with his time, whip-smart, and was a vegetarian
known for transcendental meditation.  He could easily talk about just about any
subject, ranging from sports to politics to music.  He was a family man who
suffered greatly with the loss of both a parent and a child to suicide.  Jeff
died unexpectedly on September 29, 2017 at age 68 after suffering a heart
attack and collapsing at work in front of his office door.

This paper's authors knew Jeff well.  Sergey Kitaev is Jeff's friend and his
most prolific combinatorics collaborator with 21 publications.  Anthony Mendes
is Jeff's PhD student, collaborator, and coauthor of Jeff's only book.  We are
thankful to have this opportunity to share some of Jeff's best work in
combinatorics.

\section{Symmetric functions}

Jeff's first results in combinatorics involved symmetric functions and
tableaux.  A common theme among these papers was the interpretation of the
coefficient of one symmetric function in another symmetric function as a signed
sum of combinatorial objects \cite{MR1394942}.  Jeff was then likely to
leverage this understanding to prove new results.

For instance, Jeff provided a particularly nice combinatorial interpretation
for the entries in the inverse Kostka matrix \cite{MR1034417, MR1093198}.  A
special rim hook is a sequence of connected cells in the Young diagram of an
integer partition (following Jeff's lead, we use the French convention when
drawing Young diagrams) that begins in the top left cell and travels along the
northeast edge such that its removal leaves the Young diagram of a smaller
integer partition.  A special rim hook tabloid of shape $\lambda$ and content
$\mu = (\mu_1,\dots,\mu_\ell)$ is a filling of the cells of the Young diagram of
$\lambda$ with successive special rim hooks with lengths $\mu_1,\dots,\mu_\ell$ in some order;
for example, two special rim hook tabloids of shape $(5,5,4,3,1)$ and content
$(6,6,4,2)$ are:
\begin{center}
  \hfill
  \begin{tikzpicture}[scale = .5]
    \begin{scope}
      \clip (0,0) -| (5,2) -| (4,3) -| (3,4) -| (1,5) -- (0,5) -- cycle;
      \draw [color = black!10] (0,0) grid (5,5);
    \end{scope}
    \draw [thick] (0,0) -| (5,2) -| (4,3) -| (3,4) -| (1,5) -- (0,5) -- cycle;

    \draw [teal, very thick, rounded corners] (0.5,4.5) |- ++(2,-1) |- ++(1,-1);
    \draw [color=teal, fill=teal, thick] (0.5,4.5) circle (.5ex);
    \draw [color=teal, fill=teal, thick] (3.5,2.5) circle (.5ex);

    \draw [teal, very thick] (0.5,2.5) -- (1.5,2.5);
    \draw [color=teal, fill=teal, thick] (0.5,2.5) circle (.5ex);
    \draw [color=teal, fill=teal, thick] (1.5,2.5) circle (.5ex);

    \draw [teal, very thick, rounded corners] (0.5,1.5) -| ++(4,-1);
    \draw [color=teal, fill=teal, thick] (0.5,1.5) circle (.5ex);
    \draw [color=teal, fill=teal, thick] (4.5,0.5) circle (.5ex);

    \draw [teal, very thick] (0.5,0.5) -- (3.5,0.5);
    \draw [color=teal, fill=teal, thick] (0.5,0.5) circle (.5ex);
    \draw [color=teal, fill=teal, thick] (3.5,0.5) circle (.5ex);
  \end{tikzpicture}
  \hfill
  \begin{tikzpicture}[scale = .5]
    \begin{scope}
      \clip (0,0) -| (5,2) -| (4,3) -| (3,4) -| (1,5) -- (0,5) -- cycle;
      \draw [color = black!10] (0,0) grid (5,5);
    \end{scope}
    \draw [thick] (0,0) -| (5,2) -| (4,3) -| (3,4) -| (1,5) -- (0,5) -- cycle;

    \draw [teal, very thick, rounded corners] (0.5,4.5) |- ++(2,-1);
    \draw [color=teal, fill=teal, thick] (0.5,4.5) circle (.5ex);
    \draw [color=teal, fill=teal, thick] (2.5,3.5) circle (.5ex);

    \draw [teal, very thick, rounded corners] (0.5,2.5) -- ++(3,0) |- ++(1,-1);
    \draw [color=teal, fill=teal, thick] (0.5,2.5) circle (.5ex);
    \draw [color=teal, fill=teal, thick] (4.5,1.5) circle (.5ex);

    \draw [teal, very thick, rounded corners] (0.5,1.5) -- ++(2,0) |- ++(2,-1);
    \draw [color=teal, fill=teal, thick] (0.5,1.5) circle (.5ex);
    \draw [color=teal, fill=teal, thick] (4.5,0.5) circle (.5ex);

    \draw [teal, very thick, rounded corners] (0.5,0.5) -- ++(1,0);
    \draw [color=teal, fill=teal, thick] (0.5,0.5) circle (.5ex);
    \draw [color=teal, fill=teal, thick] (1.5,0.5) circle (.5ex);
  \end{tikzpicture}
  \hfill
  \mbox{}
\end{center}
Jeff showed that the coefficient of the Schur symmetric function $s_\lambda$ in the
monomial symmetric function $m_\mu$ (those unfamiliar with these definitions are
referred to the graduate level textbook for which Jeff was a coauthor
\cite{MR3410908}) is equal to
\begin{equation*}
  \sum_{\text{special rim hook tabloids $T$ of shape $\lambda$ and content $\mu$}}
  (-1)^{\text{the number of vertical steps in $T$}}
\end{equation*}
where a vertical step is any place where a special rim hook travels down a row.

In a similar vein, let $B_{\lambda,\mu}$ be the set of all possible Young diagrams of
$\mu \vdash n$ where the rows of $\mu$ are partitioned into ``bricks'' of lengths giving
$\lambda \vdash n$.  For example, the four $T \in B_{\lambda,\mu}$ when
$\lambda = (4,2,2,1,1)$ and $\mu = (5,3,2)$ are
\begin{center}
  \hfill
  \begin{tikzpicture}[scale = .5]
    \coordinate (a) at (0,0);
    \coordinate (b) at (1,0);
    \coordinate (c) at (0,1);
    \coordinate (d) at (1,1);
    \coordinate (e) at (0,2);
    \draw [thick, rounded corners] (a) rectangle ++(1,1);
    \begin{scope}
      \clip (b) rectangle ++(3,1);
      \draw [color = black!10] (b) grid ++(4,1);
    \end{scope}
    \draw [thick, rounded corners] (b) rectangle ++(4,1);
    \draw [thick, rounded corners] (c) rectangle ++(1,1);
    \begin{scope}
      \clip (d) rectangle ++(2,1);
      \draw [color = black!10] (d) grid ++(2,1);
    \end{scope}
    \draw [thick, rounded corners] (d) rectangle ++(2,1);
    \begin{scope}
      \clip (e) rectangle ++(2,1);
      \draw [color = black!10] (e) grid ++(2,1);
    \end{scope}
    \draw [thick, rounded corners] (e) rectangle ++(2,1);
  \end{tikzpicture}
  \hfill
  \begin{tikzpicture}[scale = .5]
    \coordinate (a) at (0,0);
    \coordinate (b) at (1,0);
    \coordinate (c) at (2,1);
    \coordinate (d) at (0,1);
    \coordinate (e) at (0,2);
    \draw [thick, rounded corners] (a) rectangle ++(1,1);
    \begin{scope}
      \clip (b) rectangle ++(3,1);
      \draw [color = black!10] (b) grid ++(4,1);
    \end{scope}
    \draw [thick, rounded corners] (b) rectangle ++(4,1);
    \draw [thick, rounded corners] (c) rectangle ++(1,1);
    \begin{scope}
      \clip (d) rectangle ++(2,1);
      \draw [color = black!10] (d) grid ++(2,1);
    \end{scope}
    \draw [thick, rounded corners] (d) rectangle ++(2,1);
    \begin{scope}
      \clip (e) rectangle ++(2,1);
      \draw [color = black!10] (e) grid ++(2,1);
    \end{scope}
    \draw [thick, rounded corners] (e) rectangle ++(2,1);
  \end{tikzpicture}
  \hfill
  \begin{tikzpicture}[scale = .5]
    \coordinate (a) at (4,0);
    \coordinate (b) at (0,0);
    \coordinate (c) at (0,1);
    \coordinate (d) at (1,1);
    \coordinate (e) at (0,2);
    \draw [thick, rounded corners] (a) rectangle ++(1,1);
    \begin{scope}
      \clip (b) rectangle ++(3,1);
      \draw [color = black!10] (b) grid ++(4,1);
    \end{scope}
    \draw [thick, rounded corners] (b) rectangle ++(4,1);
    \draw [thick, rounded corners] (c) rectangle ++(1,1);
    \begin{scope}
      \clip (d) rectangle ++(2,1);
      \draw [color = black!10] (d) grid ++(2,1);
    \end{scope}
    \draw [thick, rounded corners] (d) rectangle ++(2,1);
    \begin{scope}
      \clip (e) rectangle ++(2,1);
      \draw [color = black!10] (e) grid ++(2,1);
    \end{scope}
    \draw [thick, rounded corners] (e) rectangle ++(2,1);
  \end{tikzpicture}
  \hfill
  \begin{tikzpicture}[scale = .5]
    \coordinate (a) at (4,0);
    \coordinate (b) at (0,0);
    \coordinate (c) at (2,1);
    \coordinate (d) at (0,1);
    \coordinate (e) at (0,2);
    \draw [thick, rounded corners] (a) rectangle ++(1,1);
    \begin{scope}
      \clip (b) rectangle ++(3,1);
      \draw [color = black!10] (b) grid ++(4,1);
    \end{scope}
    \draw [thick, rounded corners] (b) rectangle ++(4,1);
    \draw [thick, rounded corners] (c) rectangle ++(1,1);
    \begin{scope}
      \clip (d) rectangle ++(2,1);
      \draw [color = black!10] (d) grid ++(2,1);
    \end{scope}
    \draw [thick, rounded corners] (d) rectangle ++(2,1);
    \begin{scope}
      \clip (e) rectangle ++(2,1);
      \draw [color = black!10] (e) grid ++(2,1);
    \end{scope}
    \draw [thick, rounded corners] (e) rectangle ++(2,1);
  \end{tikzpicture}
  \hfill
  \mbox{}
\end{center}
Jeff showed that
\begin{equation}
  \label{eq:bricks}
  h_\mu = \sum_{\lambda} (-1)^{n - \ell(\lambda)} \left | B_{\lambda,\mu} \right | e_\lambda
\end{equation}
where $h_\mu$ is the homogeneous symmetric function and $e_\lambda$ the elementary
symmetric function \cite{MR1137989}.

Jeff was able to use the combinatorics of these brick tabloids to find
generating functions for permutation statistics and other objects.  His
strategy roughly followed these steps:
\begin{enumerate}
\item Define a ring homomorphism $\varphi$ on the ring of symmetric functions by
  defining $\varphi(e_{(n)})$ for all $n \geq 1$.  Since
  $e_{(1)}, e_{(2)}, \dots$ generate the ring of symmetric functions, $\varphi$ extends
  to all other symmetric functions.
\item Apply $\varphi$ to \eqref{eq:bricks} and use a sign reversing involution on the
  combinatorial objects built using brick tabloids to cancel the negative
  signs, leaving only positive fixed points.  If $\varphi$ is cleverly defined, these
  fixed points will be interesting for some reason.
\item Apply $\varphi$ to the identity
  \begin{equation*}
    \sum_{n} h_{(n)} z^n = \frac{1}{\sum_{n} (-1)^n e_{(n)} z^n}
  \end{equation*}
  to find a generating function for the fixed points.
\end{enumerate}
As an example of this idea, defining
$\varphi(e_{(n)}) = \dfrac{(-1)^n (x-1)^{n-1}}{n!}$ gives
\begin{equation*}
\sum_{n} \frac{z^n}{n!} \sum_{\sigma \in S_n} x^{\des \sigma} = \frac{x-1}{x-e^{z(x-1)}}
\end{equation*}
where $\des \sigma$ is the number of descents in the permutation $\sigma$.  As another
example, defining
$\begin{displaystyle} \varphi(e_{(n)}) = (-1)^{n-1} q^{\binom{n}{2}} \qbinom{k}{n}_q
  (x-1)^{n-1} \end{displaystyle}$ gives
\begin{equation*}
  \sum_{n} z^n \sum_{w \in \{ 0,\dots, k-1\}^*_n} x^{\des w} q^{\text{sum} w}
  = \frac{x-1}{x - \left( z - z x ; q \right)_k}
\end{equation*}
where $\{ 0,\dots, k-1\}^*_n$ is the set of words of length $n$ with letters in
$0,\dots,k-1$, $\text{sum} w$ is the sum of the integers in $w$, and we are using
the usual notations for $q$-analogues.  Other examples of this strategy can
find generating functions for linear recurrences, objects counted by the
exponential formula (permutations with restricted cycle structure, set
partitions, etc.), orthogonal polynomials such as the Chebyshev and Hermite
polynomials, and much more.  The numerous papers in this development are
thoroughly recounted in \cite{MR3410908}.

Jeff favored simple proofs and did not enjoy producing results that build on
theory or require a significant amount of mathematical overhead.  He
particularly relished proofs by bijection and sign reversing involution.  One
of Jeff's favorite proofs by bijection showed the equivalence of the definition
of the Schur symmetric functions in terms of a quotient of Vandermonde-like
determinants and the definition of the Schur symmetric function in terms of
column strict tableaux \cite{MR1357775}.  He also enjoyed his newer proof of
the Murnaghan-Nakayama rule \cite{MR3091056}.  Other interesting combinatorial
arguments are found in \cite{MR947757, MR647727, MR769977, MR694468, MR768997}.

Jeff liked computing the Littlewood-Richardson coefficients (when viewed as the
coefficient of the Schur function $s_\lambda$ in the skew-Schur function
$s_{\alpha/\beta}$) by drawing trees of standard tableaux \cite{MR769977, MR1085775,
  MR1661373}.  He used his interpretation in calculating special cases of
Kronecker coefficients (which give the number of copies of an irreducible
representation in the tensor product of two irreducible representations of a
symmetric group), a difficult problem for which formulas are only known in
certain edge cases \cite{MR1967315, MR1315363, MR1384652}.  Most of Jeff's work
here involved Schur functions $s_\lambda$ when $\lambda$ has the shape of a hook.  These
Schur functions of hook shapes were also studied in conjunction with
permutation statistics and other topics \cite{MR915951, MR1080688, MR777704,
  MR977863, MR1158791, MR1607961}.

Jeff had a good number of publications that used the plethysm of symmetric
functions and $\lambda$-ring notation \cite{MR1365453, MR1627327, MR1661367,
  MR777698, MR2035305, MR2765321}.  This is a somewhat esoteric topic that can
help when understanding the relationship between symmetric functions and the
representation theory of the symmetric group.  For an example of one such
publication, Jeff used plethystic notation to find analogues of the
Murnaghan-Nakayama rule and the calculation of Kronecker products for wreath
product groups of the form $G \wr S_n$ for a finite group $G$ \cite{MR2097322}.

\section{Enumerative combinatorics}

One of Jeff's first and most well known enumerative combinatorics results has
come to be known as Remmel's bijection machine \cite{MR676746}.  Let
$A_1,A_2,\dots$ and $B_1,B_2,\dots$ be multisets of integers such that
\begin{equation*}
  \sum \Big (\bigcup_{i \in S} A_i \Big ) = \sum \Big ( \bigcup_{i \in S} B_i \Big )
\end{equation*}
for all finite subsets $S$ of the positive integers where the union denotes a
multiset union and the sum denotes a multiset sum.  Jeff leveraged the
Garsia-Milne involution principle to find a bijection that proves the integer
partition identity
\begin{equation*}
  \left | \{ \text{$\lambda \vdash n$ with no $A_i$ in the parts} \} \right |
  = \left | \{ \text{$\lambda \vdash n$ with no $B_i$ in the parts} \} \right | .
\end{equation*}
For example, if $A_i = \{2i\}$ and $B_i = \{i,i\}$, then Remmel's bijection
machine produces a bijective proof of the identity
\begin{equation*}
  \left | \{ \text{$\lambda \vdash n$ with no even parts} \} \right |
  =   \left |\{ \text{$\lambda \vdash n$ with no repeated parts} \} \right |.
\end{equation*}
Jeff liked to point out that this proves an uncountable number of integer
partition identities bijectively and has said ``You could sit down with your
friends over a drink and say, `Hey, want to see me come up with some partition
theorems?'\,''.

Throughout his entire career, Jeff enjoyed finding $q$-analogues for
identities, having once said during a combinatorics lecture ``Let me prove one
more theorem before I $q$-analog everything in sight''.  He provided
$q$-analogues for Lagrange inversion \cite{MR873645}, sequences related to the
Fibonacci sequence \cite{MR3426219}, and even bases for the ring of symmetric
functions \cite{MR1416016}.

As an example of one such $q$-analogue, let $d_n$ denote the number of
derangements of $n$ (permutations $\sigma \in S_n$ without fixed points).  Arrange the
cycles in a permutation $\sigma \in S_n$ such that the second smallest element in each
cycle is the rightmost element in the cycle and such that cycles are ordered in
increasing order according to these second smallest elements.  Define
$\overline{\sigma}$ to be the permutation in one line notation created by removing
the parentheses in the cycles of $\sigma$.  Then if
\begin{equation*}
  d_{n,q} = \sum_{\text{derangements $\sigma \in S_n$}} q^{\inv \overline{\sigma}},
\end{equation*}
Jeff showed that
\begin{equation*}
  d_{n+1,q} = q [n]_q d_{n,q} + [n]_q d_{n-1,q} \qquad \text{and} \qquad
  d_{n+1,q} = [n+1]_q d_n + (-1)^{n+1}
\end{equation*}
for $n \geq 2$, providing $q$-analogues for well known recursions for $d_n$
\cite{MR2529717, MR576766}.

Jeff was also able to provide a new proof and a $q$-analogue of the fact that
there are $n^{n-2}$ trees on $n$ labeled vertices, a result known as Cayley's
formula.  There are numerous wonderful proofs of this theorem, including proofs
using Pr\"ufer sequences, Kirkoff's matrix tree theorem, and many more.  Jeff's
proof surpasses many of these proofs in terms of beauty and simplicity,
provides a $q$-analogue that keeps track of rises and falls in graph edges, and
can be adapted to provide an algorithm for ranking and unranking trees
\cite{MR843460, MR1382322, MR2081166}.  The proof is even short enough to
include here.

\begin{proof}
  Each of the $n^{n-2}$ functions $f : \{ 2,\dots,n-1 \} \to \{1,\dots,n\}$ can be
  represented as a directed graph on vertices $1,\dots,n$ by drawing an edge from
  $i$ to $f(i)$ for all $i$.  Draw the graph such that vertices in cycles are
  colinear with the least element in each cycle listed first and such that
  cycles are listed in decreasing order according to minimum element.  Draw any
  vertices not contained in a cycle below this line.  For example, one such
  graph is
  \begin{center}
    \begin{tikzpicture}
      \node [draw, circle, inner sep = .25ex] (12) at (0,0) {\small 12};
      \node [draw, circle, inner sep = .25ex] (4) at (2,0) {\small 4};
      \node [draw, circle, inner sep = .25ex] (7) at (3,0) {\small 7};
      \node [draw, circle, inner sep = .25ex] (8) at (4,0) {\small 8};
      \node [draw, circle, inner sep = .25ex] (3) at (6,0) {\small 3};
      \node [draw, circle, inner sep = .25ex] (1) at (8,0) {\small 1};
      \node [draw, circle, inner sep = .25ex] (9) at (0,-1) {\small 9};
      \node [draw, circle, inner sep = .25ex] (10) at (3,-1) {\small 10};
      \node [draw, circle, inner sep = .25ex] (6) at (5.25,-1) {\small 6};
      \node [draw, circle, inner sep = .25ex] (2) at (6.75,-1) {\small 2};
      \node [draw, circle, inner sep = .25ex] (5) at (8,-1) {\small 5};
      \node [draw, circle, inner sep = .25ex] (11) at (3,-2) {\small 11};
      \draw [->, ultra thick, teal] (9) -- (12);
      \draw [->, ultra thick, teal] (4) -- (7);
      \draw [->, ultra thick, teal] (7) -- (8);
      \draw [->, ultra thick, teal] (8) .. controls ++(-0.5,.75) and ++(.5,.75) .. (4);
      \draw [->, ultra thick, teal] (10) -- (7);
      \draw [->, ultra thick, teal] (11) -- (10);
      \draw [->, ultra thick, teal] (3) .. controls ++(.5,.75) and ++(-.5,.75) .. (3);
      \draw [->, ultra thick, teal] (6) -- (3);
      \draw [->, ultra thick, teal] (2) -- (3);
      \draw [->, ultra thick, teal] (5) -- (1);
    \end{tikzpicture}
  \end{center}

  Bijectively change this graph into a tree by connecting the cycles from left
  to right, erasing loops, and undirecting edges.  Doing this to the above
  gives
  \begin{center}
    \begin{tikzpicture}
      \node [draw, circle, inner sep = .25ex] (12) at (0,0) {\small 12};
      \node [draw, circle, inner sep = .25ex] (4) at (2,0) {\small 4};
      \node [draw, circle, inner sep = .25ex] (7) at (3,0) {\small 7};
      \node [draw, circle, inner sep = .25ex] (8) at (4,0) {\small 8};
      \node [draw, circle, inner sep = .25ex] (3) at (6,0) {\small 3};
      \node [draw, circle, inner sep = .25ex] (1) at (8,0) {\small 1};
      \node [draw, circle, inner sep = .25ex] (9) at (0,-1) {\small 9};
      \node [draw, circle, inner sep = .25ex] (10) at (3,-1) {\small 10};
      \node [draw, circle, inner sep = .25ex] (6) at (5.25,-1) {\small 6};
      \node [draw, circle, inner sep = .25ex] (2) at (6.75,-1) {\small 2};
      \node [draw, circle, inner sep = .25ex] (5) at (8,-1) {\small 5};
      \node [draw, circle, inner sep = .25ex] (11) at (3,-2) {\small 11};
      \draw [ultra thick, teal] (9) -- (12);
      \draw [ultra thick, teal] (4) -- (7);
      \draw [ultra thick, teal] (7) -- (8);
      \draw [ultra thick, teal] (10) -- (7);
      \draw [ultra thick, teal] (11) -- (10);
      \draw [ultra thick, teal] (6) -- (3);
      \draw [ultra thick, teal] (2) -- (3);
      \draw [ultra thick, teal] (5) -- (1);
      \draw [ultra thick, teal] (12) -- (4);
      \draw [ultra thick, teal] (8) -- (3);
      \draw [ultra thick, teal] (3) -- (1);
    \end{tikzpicture}
  \end{center}
  The directed graph was carefully drawn in the prescribed manner as to make
  this process bijective; the details are left to the reader or see
  \cite{MR843460}.
\end{proof}

Another one of Jeff's pet topics was rook theory.  A rook board $B$ of size $n$
is a sequence of columns of cells of heights $(0,\dots,n-1)$ atop columns of
infinite depth that contain $n$ non-attacking rooks.  Here, similar to but not
exactly the same as chess, rooks attack all cells below and to the right.  For
instance, one board $B$ of size $4$ is
\begin{center}
  \begin{tikzpicture}[scale = .5]
    \begin{scope}
      \clip (0,0) -| (1,1) -| (2,2) -| (3,3) -- (4,3) -- (4,0) -- cycle;
      \draw [color = black!10] (0,0) grid (4,4);
    \end{scope}
    \draw [thick] (0,0) -| (1,1) -| (2,2) -| (3,3) -- (4,3) -- (4,0) -- cycle;
    \draw [color = black!10] (0,0) grid (4,-6.75);
    \draw [thick] (0,0) -- (4,0);
    \draw [thick] (0,0) -- (0,-6.75);
    \draw [thick] (4,0) -- (4,-6.75);
    \node at (-0.5,-0.5) {\small{1}};
    \node at (-0.5,-1.5) {\small{2}};
    \node at (-0.5,-2.5) {\small{3}};
    \node at (-0.5,-3.5) {\small{4}};
    \node at (-0.5,-4.5) {\small{5}};
    \node at (-0.5,-5.5) {\small{6}};
    \node at (-0.5,-6.25) {\small{$\vdots$}};
    \node (r1) at (2.5,1.5) {\rook};
    \node (r2) at (1.5,0.5) {\rook};
    \node (r3) at (0.5,-2.5) {\rook};
    \node (r4) at (3.5,-4.5) {\rook};
    \draw [color = teal, line width = .5ex] (r1) -- ++(1.5,0);
    \draw [color = teal, line width = .5ex] (r1) -- ++(0,-8.25);
    \draw [color = teal, line width = .5ex] (r2) -- ++(2.5,0);
    \draw [color = teal, line width = .5ex] (r2) -- ++(0,-7.25);
    \draw [color = teal, line width = .5ex] (r3) -- ++(3.5,0);
    \draw [color = teal, line width = .5ex] (r3) -- ++(0,-4.25);
    \draw [color = teal, line width = .5ex] (r4) -- ++(0,-2.25);
  \end{tikzpicture}
\end{center}
Let $\inv B$ be the number of non-attacked cells and $\max B$ is the row with
the lowest rook (the above example has $\inv B = 6$ and $\max B = 5$).  One of
Jeff's first papers on rook theory showed that if
\begin{equation*}
  S_{n,k} (q)
  = \sum_{\text{placements of $n-k$ rooks on a board $B$ of shape $(0,\dots,n-1)$}} q^{\text{inv} B},
\end{equation*}
then $S_{n,k}(q)$ satisfies
\( S_{n-1,k}(q) = q^{k-1} S_{n,k-1}(q) + [k]_q S_{n,k}(q) \) and
\begin{equation*}
  \sum_k \frac{S_{n,k}(q) \; [k]_q! \; t^k}{(1-t q^0) \cdots (1-t q^k)} =
  \sum_{B} t^{\text{max} B} q^{\text{inv B}}
  = \sum_{\sigma \in S_n} \frac{ q^{\text{maj} \sigma} t^{\text{des} \sigma+ 1}}{(1-t q^1)\cdots(1-t q^n)}
\end{equation*}
where $\maj{\sigma}$ and $\des{\sigma}$ are the major index and descent statistics for
permutations in the symmetric group $S_n$.  Jeff and coauthors generalized this
in a myriad of ways, finding permutation enumeration results for groups other
than $S_n$ and results and identities for various shapes of boards
\cite{MR3176194, MR3508003, MR2244138, MR834272, MR1868975, MR3930459,
  MR2511012, MR2426148, MR2114188}.

Jeff also dabbled with other enumerative combinatorics topics, including
perfect matchings \cite{MR2143200, MR2025077} and shuffles \cite{MR951014,
  MR3678665, MR4099466}.

\section{Patterns in combinatorial structures}

Studying the appearance of patterns in combinatorial structures (primarily
permutations and words) was a significant part of Jeff's research during the
final 10 years of his life.  Jeff authored 60 papers on the subject (out of his
115 articles published since 2007) and he wrote the foreword for the only
comprehensive book on the subject \cite{MR3012380}.

An occurrence of a pattern $\tau$ in a permutation $\sigma$ is ``classically'' defined
as a subsequence of $\sigma$ whose elements are in the same relative order as those
in $\tau$, and Jeff published a couple of papers about such patterns
\cite{MR3382607, MR4028903}.  However, the notion of a pattern has been
extended to other settings many times in the literature, and Jeff was behind
several of innovations (for example, \cite{MR2240770, MR2601315, MR3325551,
  MR3338849}).

The notion of a quadrant marked mesh pattern, introduced by Jeff in
\cite{MR2914898}, resulted in a large program of research by Jeff and coauthors
in a series of papers \cite{MR3118905, MR3064041, MR3119706, MR3345153,
  MR3418503, MR3843149}.  Let $\sigma = \sigma_1 \cdots \sigma_n$ be a permutation in the symmetric
group $S_n$ written in one-line notation.  Then $\sigma_i$ matches the quadrant marked mesh pattern
$\text{MMP}(a,b,c,d)$ if the number of points $(j,\sigma_j)$ in the four quadrants
with origin at $(i,\sigma_i)$ satisfies the inequalities as depicted below:
\begin{center}
\begin{tikzpicture}[scale = .25]
  \draw [color = black!10] (1,1) grid (7,7);
  \draw [color = black!50, thick] (4,1) -- (4,7);
  \draw [color = black!50, thick] (1,4) -- (7,4);
  \draw [color=black, fill=black, thick] (4, 4) circle (1ex);
  \node at (5.5,5.5) {\color{black} {$\geq a$}};
  \node at (5.5,2.5) {\color{black} {$\geq d$}};
  \node at (2.5,2.5) {\color{black} {$\geq c$}};
  \node at (2.5,5.5) {\color{black} {$\geq b$}};
\end{tikzpicture}
\end{center}
For example, the `6' in {4\,7\,1\,5\,6\,9\,2\,8\,3} satisfies $MMP(2,0,3,1)$:
\begin{center}
\begin{tikzpicture}[scale = .33]
  \draw [color = black!10] (1,1) grid (9,9);

  \draw [color = black!50, thick] (5,1) -- (5,9);
  \draw [color = black!50, thick] (1,6) -- (9,6);

  \draw [color=teal, fill=teal, thick] (1, 4) circle (1ex);
  \draw [color=teal, fill=teal, thick] (2, 7) circle (1ex);
  \draw [color=teal, fill=teal, thick] (3, 1) circle (1ex);
  \draw [color=teal, fill=teal, thick] (4, 5) circle (1ex);
  \draw [color=teal, fill=teal, thick] (5, 6) circle (1ex);
  \draw [color=teal, fill=teal, thick] (6, 9) circle (1ex);
  \draw [color=teal, fill=teal, thick] (7, 2) circle (1ex);
  \draw [color=teal, fill=teal, thick] (8, 8) circle (1ex);
  \draw [color=teal, fill=teal, thick] (9, 3) circle (1ex);

  \node at (1,0.25) {\small{1}};
  \node at (2,0.25) {\small{2}};
  \node at (3,0.25) {\small{3}};
  \node at (4,0.25) {\small{4}};
  \node at (5,0.25) {\small{5}};
  \node at (6,0.25) {\small{6}};
  \node at (7,0.25) {\small{7}};
  \node at (8,0.25) {\small{8}};
  \node at (9,0.25) {\small{9}};

  \node at (0.25,1) {\small{1}};
  \node at (0.25,2) {\small{2}};
  \node at (0.25,3) {\small{3}};
  \node at (0.25,4) {\small{4}};
  \node at (0.25,5) {\small{5}};
  \node at (0.25,6) {\small{6}};
  \node at (0.25,7) {\small{7}};
  \node at (0.25,8) {\small{8}};
  \node at (0.25,9) {\small{9}};
\end{tikzpicture}
\end{center}
As a sample of one of his results, Jeff showed that
\begin{equation*}
  1 + \sum_{ n = 1}^\infty \frac{t^{2n}}{(2n!)}
  \sum_{\text{$\sigma \in S_{2n}$ is alternating}}
  q^{\text{mmp}^{(1,0,0,0)}(\sigma)} = \sec(q t)^{1/q},
\end{equation*}
where $\text{mmp}^{(1,0,0,0)}(\sigma)$ is the number of elements in $\sigma$ that match
$\text{MMP}(1,0,0,0)$, thereby refining Andr\'e's classical result on
alternating permutations.  A similar generating function for permutations of
odd length is $\int_0^t (\sec(q z))^{1+\frac{1}{q}} \, dz$.  These elegant results
show yet again that Jeff could $q$-analgue almost anything!

A major stream of Jeff's research on patterns was related to consecutive
patterns \cite{MR3426220, MR3266909, MR3106443, MR2601804, MR2266893,
  MR3903152, MR3262239}, where Jeff leveraged his knowledge of symmetric
functions.  An occurrence of a consecutive pattern is always formed by
contiguous elements in a permutation or word, but otherwise it is an occurrence
of a classical pattern. A particular result in this direction from
\cite{MR3106443} is showing that
\begin{equation*}
  \sum_{n\geq 0}\frac{t^n}{n!}\sum_{\sigma \in S_n}x^{\text{LRmin}(\sigma)}y^{1+\des \sigma}
  =\left(1+\sum_{n = 1}^\infty U_{\tau,n}(y)\frac{t^n}{n!}\right)^{-x}
\end{equation*}
where $\sigma$ runs over all permutations of length $n$ which avoid a consecutive
pattern $\tau$ having only one descent and the element 1 in the first position,
LRmin is the left-to-right minima statistic, and the coefficients
$U_{\tau,n}(y)$ satisfy simple recursions.

Jeff studied the bivincular pattern \patterna related to the interval orders
and ascent sequences encoding them, as well as to several other remarkable
combinatorial objects \cite{MR2847912, MR2832334, MR3668878}.  An occurrence of
\patterna in a permutation is an occurrence of the pattern $2\,3\,1$ in which
the first and second elements are next to each other, and the first element is
one more than the last element.

In \cite{MR2847912}, Jeff showed that the ordinary generating function for the
number of \patterna-avoiding permutations with at most $k$ consecutive elements
in decreasing order that are next to each other in value (of the form
$a(a-1)(a-2)\cdots$) is given by
\begin{equation*}
  \sum_{n = 0}^\infty \prod_{i=1}^{n} \left(1-\left(\frac{1-t}{1-t^k}\right)^i\right).
\end{equation*}

In \cite{MR2832334}, Jeff proposed an interesting conjecture that the ordinary
generating function for the number of \patterna-avoiding permutations with the
leftmost decreasing run of size $k$ (controlled by the variable $z$) is
\begin{equation*}
  \sum_{n = 0}^\infty \prod_{i=1}^{n}(1-(1-t)^{i-1}(1-zt)).
\end{equation*}
This former conjecture refines an important enumerative result in \cite{MR2652101}.

Jeff gave an unexpected application of patterns in graph representations
\cite{MR3367296}. The basic idea is that graphs can be encoded by words where
the edge relations are determined by occurrences of a fixed pattern in a word.
This is a far-reaching generalization of the notion of a word-representable
graph \cite{wr-wiki}.  Jeff went even further, and communicated (less than 3
months prior to his death) the idea of tolerance to occurrences of patterns
defining edges/non-edges in graph representations. This idea was implemented in
\cite{MR3960511} where it was shown that every graph is 2-11-representable
(leaving open the challenging question whether every graph is
1-11-representable).

Another topic worth mention is the notion of a generalized factor order
on words \cite{MR2721539,MR2792158}.

Jeff's extensive work on patterns in combinatorial objects is only touched upon
here, although we have pointed the reader to many references throughout this
paper of his best work.  Interested parties are certainly encouraged to read
some of these papers.

\bibliography{Remmel}
\bibliographystyle{plain}

\end{document}